\documentclass[9pt]{article}
\usepackage{amsmath}
\usepackage[utf8]{inputenc}
\usepackage [autostyle, english = american]{csquotes}
\usepackage{enumitem}
\usepackage{authblk}
\usepackage{lineno}
\usepackage[toc,page]{appendix}
\usepackage{amsthm}
\usepackage{changepage}
\usepackage{aliascnt}
\usepackage{theoremref,amsthm}
\usepackage{morefloats}
\usepackage{amsfonts}
\usepackage{amssymb}
\usepackage{graphics}
\usepackage{graphicx}
\usepackage{listings}
\usepackage{xcolor}
\usepackage{float}
\usepackage{verbatim}
\usepackage{geometry}
\usepackage{centernot}
\usepackage[english]{babel}
\usepackage{setspace}
\usepackage{thmtools}
\usepackage{xcolor}
\usepackage[bottom]{footmisc}
\usepackage{euscript}
\usepackage{mathtools}
\DeclareMathAlphabet{\itbf}{OML}{cmm}{b}{it}

\newcommand{\bea}{\begin{eqnarray*}}
\newcommand{\eea}{\end{eqnarray*}}
\newcommand{\bean}{\begin{eqnarray}}
\newcommand{\eean}{\end{eqnarray}}
\newcommand{\p}{\partial}
\newcommand{\f}{\frac}
\newcommand{\s}{\sqrt}
\newcommand{\ds}{\displaystyle}
\newcommand{\no}{\nonumber}
\newcommand\ov{\overline}
\newcommand{\ri}{\rightarrow}
\newcommand{\vp}{\varphi}

\newtheorem{lem}{Lemma}[section]
\newtheorem{prop}{Proposition}[section]

\newtheorem{thm}{Theorem}[section]

\newcommand{\RR}{\mathbb{R}}

\newcommand{\NN}{\mathbb{N}}
\newcommand{\CC}{\mathbb{C}}

\renewcommand{\ov}[1]{\overline{#1}}

\newcommand{\by}{{\itbf y}}

 \author{ Darko  Volkov 
\thanks{Department of Mathematical Sciences,
Worcester Polytechnic Institute, Worcester, MA 01609. Corresponding author email: darko@wpi.edu. 
} }

\begin{document}
		\title{Optimal decay rates 
		in Sobolev norms for singular values
	of 	integral operators }
		\maketitle

\begin{abstract}
The regularity of integration kernels forces decay rates of singular values
of associated integral operators. 
This is well-known for symmetric operators with kernels defined
on $(a,b) \times (a,b)$, where $(a,b)$ is an interval.
Over time, many authors
  have  studied this case in detail  
\cite{GohbergKrein, ha1986, chang1999eigenvalues, 
little1984eigenvalues, hille1931characteristic}. The case of spheres
has also been resolved  \cite{castro2012eigenvalue}.
A few authors have examined the  higher dimensional case or 
the case of
manifolds \cite{birman1977estimates, delgado2021schatten}.
Typically, these authors have provided decay estimates 
of  singular values in $l^p$ norms, $p \geq 1$ or in case
of faster decay due to regularity, $l^p$ quasi-norms, $0 < p < 1$. 
With that approach, it is straightforward to show that their estimates
are optimal using periodic kernels obtained from Fourier series.
Our new approach for deriving decay estimates of these singular values
uses Weyl's asymptotic formula for Neumann eigenvalues  
\cite{weyl1912asymptotische} that we combine to an appropriately defined 
inverse Laplacian. 
We obtain decay estimates in the form $n^\alpha$ where
for the $n$ -th singular value where $\alpha $ depends on  dimension and on the Sobolev regularity of the kernel.
Since we are interested in optimal estimates in case of regular kernels,
 instead  of writing 
an upper bound  by a constant times $n^\alpha$, we use the 
singular values of the kernel obtained by differentiation. 
While in  \cite{birman1977estimates, delgado2021schatten} $l^p$
estimates were proven to be optimal by simply considering  periodic 
Fourier series, these Fourier series do not provide
sharp results for our estimates. Instead, series of
Neumann eigenfunctions for the Laplacian that are specific to the domain
of interest are used to
prove that our decay estimates are optimal.
Finally, we cover the case of real analytic kernels where we 
are also able to derive optimal estimates.
\end{abstract}

\section{Introduction} \label{intro}
Let $k:  \Omega' \times \Omega \ri \RR$ be a function of
$(x,y)$. 
The integral operator
 $K$ associated to  kernel the $k$ is defined
by setting 
 $(K \phi) (x) = \int_\Omega k(x,y) \phi (y) dy$.
Here $\Omega$ and $\Omega'$ can be open subsets of $\RR^d$ and
$\RR^{d'}$ or manifolds, with or without boundaries.
If $k$ is in $L^2(\Omega' \times \Omega )$ then $K$ defines a compact operator from $L^2(\Omega)$ to $L^2(\Omega')$.
It can be observed in computations
that the more regular $k$ is, the faster the singular values of $K$ decay. 
The decay rate of the singular values of $K$ relates to the ill-posedness of the 
inverse (or generalized  inverse) of $K$. This is of great importance 
for building regularizing schemes for the numerical solutions to inverse 
problems \cite{kress1989linear, vogel2002computational}.
It is also important in the construction of sampling type algorithms
for solving inverse problems
where a regularizing term is dilated by a parameter
which is modeled as a random variable following a
hyperprior distribution
\cite{calvetti2018inverse, volkov2022parallel}. 
Having an estimate of  decay rates of singular values
 also has implications regarding the possibility of 
building neural networks for the numerical solution to 
inverse problems stemming from the inversion of integral 
operators \cite{VolkovStability2023}.
However, relating the decay rate of singular values to the regularity of the kernel is only thoroughly understood in dimension one,
 equivalently when
$d=d'=1$.
For example,  Gohberg and Krein provide a
detailed account of decay properties 
in section III.10 of their textbook  \cite{GohbergKrein}
in the case where $\Omega=\Omega'=(a,b)$.
They 
obtain an
estimate for the decay of singular values if a first derivative in mean exists. They then 
suggest how to improve this estimate if higher order derivatives exist.
The  idea 
of using the integral operator on $[0,1]$
$Jf(x)=\int_{x}^1 f(y) dy$ and its adjoint
suggested by Gohberg and Krein  was later 
used by  Ha \cite{ha1986}. Ha gives a complete proof of decay estimates
of singular values of positive definite kernels with continuous derivatives.
From there, Ha and collaborators 
\cite{chang1999eigenvalues} considered the case of weak derivatives and special boundary conditions.
In his masters thesis \cite{levine2023decay},  Levine attempted to generalize this method to the 
higher dimensional case but failed even in the case $\Omega = (0,1)^2$.
We believe that this is due to the fact  that while the subspace of functions
$f$
such that $\int_0^1 k(1,y) f(y) dy =0$ has co-dimension 1 in the case where
 $\Omega = (0,1)$, the analog  in the case where
 $\Omega = (0,1)^2$ that would be used 
for eliminating boundary terms when integrating by parts is an infinite dimensional 
 subspace.\\
In \cite{castro2012eigenvalue}, Castro and Menegatto were
able to avoid that difficulty by considering integral operators on a sphere. 
Indeed, a sphere is a manifold without boundary so there are no boundary
terms when integrating by parts. Their approach uses the surface Laplacian on spheres.
In this study we are able to adjust this idea to integral operators on 
bounded 
open sets or manifolds with a boundary. Our idea is the repeated use 
of an inverse Laplacian defined from Neumann eigenvectors.
This inverse Laplacian defines a compact operator whose eigenvalues can be estimated. In fact, Weyl was the first to systematically study the growth rate
of Neumann eigenvalues \cite{weyl1912asymptotische}.
This growth rate and the use of an appropriate inverse Laplacian proved essential in our approach. Our main estimate for the decay rate of singular 
values is stated in Theorem \ref{thm1}.
In our analysis, we found that considering classes of kernels  with continuous 
derivatives  is awkward.
It is important to consider kernel with integrable weak derivatives
to obtain optimal results.
Roughly speaking, we find that if the kernel is in $L^2(\Omega') \times 
H^p(\Omega)$, where $\Omega$  is a bounded open set of
$\RR^d$,  then the $n$-th singular value of $K$ decays
faster than $n^{-\f{p}{d}}$. 
The 
 \enquote{faster than} 
  rate can be in turn estimated 
	by  the singular values attached
to operators formed from weak derivatives of order $p$.
In particular, Theorem  \ref{thm1} implies that the 
$n$-th singular value of $K$ is  $o(n^{-\f{p}{d}})$
as shown by Ha \cite{ha1986}
if $d=1$ and the derivatives exist in the space of continuous 
functions on $[0,1]$.
 In  section \ref{proof of thm2} 
we show that the decay rate  \eqref{main est}
in Theorem \ref{thm1} of this paper
 is optimal if $\Omega $ is a domain with a
boundary of class
$C^p$. If the boundary $\p \Omega$
is only Lipschitz regular  then the estimate is almost 
optimal: for every $\alpha>0$ there is a kernel  in $L^2(\Omega) \times 
H^p(\Omega)$ such that  the $n$-th singular value of the associated
integral operator decays
more slowly  than $n^{-\f{p}{d} -\alpha}$.
Authors such as Birman and Solomyak and Delgado and Ruzhansky
\cite{birman1977estimates, delgado2021schatten}
followed a different approach: instead of directly estimating the decay rate of the sequence $s_n$ of singular values of an integral operator,
they chose to estimate the $l^p$ norm of $s_n $ and related norms.
This approach is particularly fruitful for studying kernels 
that are in the space $L^q(\Omega') \times L^r(\Omega)$.
For more regular kernels, they consider $l^p $ 
pseudo-norms of $s_n$, $0<p<1$.
Here, by more regular we mean that some weak derivatives exist.
The estimate \eqref{main est} in this paper most closely resembles estimate 5.12 in \cite{delgado2021schatten}.
For that kind of estimate, optimality can be shown by setting up
kernels derived from periodic Fourier series as pointed out in
\cite{birman1977estimates}, shortly after Corollary 4.2.
It turns out that our more direct estimate 
 \eqref{main est} is also optimal, however, optimality cannot be 
shown using periodic Fourier series. 
Instead, we have to resort to eigenfunction expansions for 
the Laplace Neumann problem specific to the domain $\Omega$.
This is the subject of Theorem \ref{thm2}.
The global regularity of these eigenfunctions depend on the 
smoothness of the boundary $\p \Omega$.
If $\p \Omega $ is only Lipschitz regular, then kernels derived from periodic Fourier series only provide quasi optimality 
as explained in the proof of Theorem \ref{thm2}.
We note that Birman and Solomyak provide in \cite{birman1977estimates}  an extensive study
of classes of integral operators and asymptotic estimates of their singular values. Their proof method is based on polynomial approximations
on cubes. In our view the English translation of their work is very dense with a difficult to follow presentation and many hard to access references. 
We believe that our method based on Weyl's estimates and an appropriate
inverse Laplacian
 offers instead a  more straightforward derivation which works
well with Sobolev norms.
To complete our literature review we also refer to
 \cite{hille1931characteristic} and \cite{birman1977estimates}.
While Hille and Tamarkin  \cite{hille1931characteristic} entirely focus on eigenvalues in the one-dimensional case, 
 they cover interested cases including 
H\"older regularity and the analytic case. 
\\
Finally, our last result concerns kernels $k(x,y)$, $(x,y) \in \Omega' \times 
\ov{\Omega}$,
  that are analytic in $y$ as a function 
valued in $L^2(\Omega')$.
We are able to generalize Little and Reade's results pertaining to 
 the one-dimensional  case ($d=d'=1$)
	\cite[Theorem 15.20]{kress1989linear}, \cite{little1984eigenvalues}
	to higher dimensions. In our multi-dimensional case we find
	that singular values of the corresponding integral operator 
	decay faster than
	$\tau^{n^{\f{1}{d}}}$ for some $\tau$ in $(0,1)$.
	Conversely, we show that for any bounded open set $\Omega$ and any
	$\tau$ in $(0,1)$ there is a kernel $k(x,y)$, analytic in $y$
	such that the  $n$-th singular value of the associated integral operator 
	is bounded below by a positive constant times $\tau^{n^{\f{1}{d}}}$.
\\
In a last paragraph, we  indicate how our results are easily generalized to the case of  
manifolds, with or without boundary, as long as they can be covered by finitely many open bounded sets where local charts can be defined.\\
Here are the three theorems that we prove in this paper
where we use the standard notation $s_n(K)$
for the $n$-th singular value of a compact operator $K$.

\begin{thm}\label{thm1}
Let $\Omega$ be a bounded open set of $\RR^d$ such that $\p \Omega$
is Lipschitz regular. Let $\Omega'$ be an open set of $\RR^{d'}$
and $k: \Omega' \times\Omega \ri \RR$
 be in $L^2(\Omega') \times H^p(\Omega)$.
Define the linear operator
\bea 
K:  L^2(\Omega) \ri L^2(\Omega') , \no \\
K \varphi(x) = \int_\Omega k(x,y) \varphi(y) dy.
\eea
Set 
	\bea
K_{i_1, ..., i_p}:  L^2(\Omega) \ri L^2(\Omega') , \\
K_{i_1, ..., i_p} \varphi(y) = \int_\Omega\p_{i_1} ... \p_{i_p}
 k(x,y)  \varphi(y) dy,
\eea
where $1\leq i_1, ..., i_p \leq d$ are integers and
all derivatives $\p_{i_1}, ... ,\p_{i_p}$ are in $y$.
Then
	\bean \label{main est}
	 s_n(K) \leq  C n ^{-\f{p}{d}} \sum_{1 \leq i_1, ..., i_p \leq d}
	s_m (K_{i_1, ..., i_p} ) , \quad  \quad m=[ \f{n}{ d^p (p+1)} ],
	\eean
where the constant $C$ does not depend on $n$ and $[\, ]$ is the integer part function. In particular
\bean \label{particular 1}
\sum_{n \geq 1 } n ^{\f{2p}{d}} s_n(K)^2 < \infty
\eean
and
\bean \label{particular 2}
 s_n(K) = o(n^{- \f{p}{d} -\f12}).
\eean
\end{thm}

\begin{thm}\label{thm2}
If $\p \Omega $ is  Lipschitz-regular 
estimate \eqref{main est}
is quasi-optimal in the sense that
then for any $\alpha >\f{1}{2}$ 
	we can find 
	 a kernel $k$  in $L^2(\Omega') \times H^p(\Omega)$, an associated operator $K$  and
a	positive constant $c$
	such that 
\bean \label{main est 2.1}
	   c n ^{-\f{p}{d} - \alpha}  \leq s_n(K) 
	\eean
	 for all positive integers $n$.\\
If $\p \Omega $ is of class $C^p$ 
estimate \eqref{main est}
is optimal in the sense that there is a kernel $k$  in $L^2(\Omega') \times H^p(\Omega)$ and an associated operator $K$ such that 
for some positive constant $c$ 
\bean \label{main est 2}
	   c n ^{-\f{p}{d}} \sum_{1 \leq i_1, ..., i_p \leq d}
	s_m (K_{i_1, ..., i_p} ) \leq s_n(K) , \quad  \quad m=[ \f{n}{ d^p (p+1)} ],
	\eean
	 for all positive integers $n$.\\
\end{thm}

\begin{thm}\label{thm3}
Let $\Omega$ be a bounded open set of $\RR^d$.
  Let $\Omega'$ be an open set of $\RR^{d'}$
and $(x,y)\ri k(x,y)$ be in $L^2(\Omega' \times\Omega)$.
Assume that $k$ is analytic in a neighborhood of $\ov{\Omega}$ as a function valued in $L^2(\Omega')$.
Then there exists $\tau $ in $(0,1)$ and $C>0$ 
such that
\bean \label{main est 3}
s_n(K) \leq C \tau^{n^{\f{1}{d}}},
\eean
for all positive integers $n$.
This estimate is optimal in the sense that  
for any open set $\Omega$ and any $\tau$ in $(0,1)$
 there is a  kernel $k$ which is analytic in a neighborhood of $\ov{\Omega}$ as a function valued in $L^2(\Omega')$,
and $c>0$ 
such that
\bean \label{main est 4}
 c  \tau^{n^{\f{1}{d}}}\leq s_n(K) ,
\eean
for all positive integers $n$.
\end{thm}

\section{Singular values of compact operators: a review of elementary properties }
In this section, all Hilbert spaces are assumed to be separable.
Let $H_1, H_2$ be two such spaces. 
Let $T: H_1 \ri H_2$ be a linear compact operator.
By definition, the sequence $s_n(T)$ of the singular values of $T$ is the sequence of square roots of  eigenvalues  
of $T^*T$ given in decreasing order, with possible repetitions to account for  multiplicity of eigenvalues. 
 It is well-known that $T^*T$ and $TT^*$ have the same non-zero eigenvalues, counted with multiplicity,
thus 
\bean \label{dual}
s_n(T)=s_n(T^*),
\eean
if $s_n(T)$ or $s_n(T^*)$ is non-zero.
\begin{prop}\label{min max prop}
Fix an integer $n\geq 1 $ and let ${\cal S}_n$ be the set of all subspaces of $H_1$ with co-dimension less or equal than   $n-1$. 
Then
\bean \label{min max}
s_n(T) = \min_{X \in {\cal S}_n} \max_{ x\in X, \| x \|=1} \| T x\|.
\eean
\end{prop}
\textbf{Proof:}
For $x \in H_1$, $\| T x\|^2 = <T^*Tx, x>$.
As  $T^*T$ is a compact, symmetric, and non-negative operator
this results from \cite[Theorem III.9.1]{gohberg2013basic}
if ${\cal S}_n$ is  the set of all subspaces of $H_1$ with co-dimension 
equal to   $n-1$. 
If $X'$ is a subspace of $H_1$ with co-dimension
$n'-1 \leq n-1$ then 
\bea
s_{n'} (T) \leq \max_{x \in X', \| x \| =1} 
\| T x  \|.
\eea
  As $ s_{n} (T) \leq s_{n'} (T)$, the result follows.
$\Box$

\begin{prop}\label{restriction}
Let $H_1, H_2$ be two  Hilbert spaces, 
$T: H_1 \ri H_2$, 
 be a  compact linear  operator, $n, q$ two positive integers, and $H_3$ a subspace
of $H_1$ of co-dimension $q$.
Then
\bean \label{restric ineq}
s_{n+q} (T) \leq s_n(T_{| H_3}).
\eean
\end{prop}
\textbf{Proof:}
Let $\{ v_1, ..., v_q \}$ be an orthonormal basis of $H_3^\perp$.
Let $\{ v_{q+1}, ..., v_{q + n} \}$ be $n-1$ orthonormal vectors in $H_3$ such that
\bea
 s_n(T_{| H_3}) =  \max_{ x\in \{ v_{q+1}, ..., v_{q + n} \}^\perp
\cap H_3, \| x \|=1} \| T x\|.
\eea
These $n-1$ vectors exist thanks to Proposition \ref{min max prop}.
Since the co-dimension of $ \{ v_{q+1}, ..., v_{q + n} \}^\perp
\cap H_3 $ in $H_1$ is $n+q -1$, the result follows by Proposition \ref{min max prop}.
$\Box$

\begin{prop} \label{norm and sing}
Let $H_1, H_2, H_3$ be three  Hilbert spaces 
Let $S: H_1 \ri H_2$, 
$T: H_2 \ri H_3$, 
 be continuous  linear  operators.
If $S$ is compact then $s_n(TS) \leq \| T \| s_n (S)$, $n \geq 1$,
if $T$ is compact then $s_n(TS) \leq \| S \| s_n (T)$, $n \geq 1$.
\end{prop}
\textbf{Proof:}
The first inequality is clear with \eqref{min max} applied to $TS$.
The second inequality  results from the duality relation \eqref{dual}.
$\Box$\\
A special case of the following proposition appears in
\cite{ha1986, castro2012eigenvalue}. We have to use a slightly more general version that applies to 
operators between Hilbert spaces.  
\begin{prop}
Let $H_1, H_2$ be  two Hilbert spaces 
Let $S: H_1 \ri H_2$, 
$T: H_2 \ri H_3$, 
 be compact   linear  operators, and $p, q \geq 1$ two integers.
Then
\bean
s_{p+q-1} (TS)\leq s_p(S)  s_q(T).  \label{pq}
\eean
\end{prop}
\textbf{Proof:}
Let $v_1, ..., v_{p-1}$ be orthonormal eigenvectors vectors in 
$H_1$ corresponding to the first $p-1$ eigenvalues of $S^*S$.
Let $\tilde{H} =\{v_1,...,v_{p-1} \}^\perp$. Note that $\tilde{H} $
is its own Hilbert space. Let $w_1, ..., w_{q-1}$ be $q-1$ independent
vectors in $\tilde{H} $ and 
$$
X=\{ x \in \tilde{H}: <x,w_1>=...=<x,w_{q-1}>=0 \}.
$$
Note that $X$ is a subspace of $H_1$ of co-dimension $p+q-2$.
By \eqref{min max},
$s_{p+q-1} (TS)\leq \max_{x \in X, \| x\|=1} 
\| TS x \|$.
Let $\tilde{S} $ be the restriction of $S$ to $\tilde{H}$.
By Proposition \ref{norm and sing},
$$
 \max_{x \in X, \| x\|=1} 
\| TS x \| =
 \max_{x \in X, \| x\|=1} 
\| T\tilde{S} x \| \leq
\| \tilde{S} \| s_q ( T ) =  s_p ( S) s_q ( T ).
$$
$\Box$
\begin{prop}
Let $H_1, H_2$ be  two Hilbert spaces 
Let $S, T: H_1 \ri H_2$, 
 be compact   linear  operators, and $p, q \geq 1$ two integers.
Then
\bean
s_{p+q-1} (S +T)\leq s_p(S) + s_q(T).  \label{pq plus}
\eean
\end{prop}
\textbf{Proof:}
Let $v_1, ..., v_{p-1}$ be orthonormal eigenvectors vectors in 
$H_1$ corresponding to the first $p-1$ eigenvalues of $S^*S$
 and $w_1, ..., w_{q-1}$  be orthonormal eigenvectors vectors in 
$H_1$ corresponding to the first $q-1$ eigenvalues of $T^*T$.
Denote $V =\{v_1, ..., v_{p-1}, w_1, ..., w_{q-1} \}^\perp$,
$V_1 =\{v_1, ..., v_{p-1} \}^\perp$,
$V_2 =\{w_1, ..., w_{q-1} \}^\perp$. Clearly, $V \subset V_1$, $V \subset V_2$.
It follows that
\bea
\max_{x \in V, \| x\| =1} \| (S+T) x  \| & \leq &
\max_{x \in V, \| x\| =1} \| S x \| +
\max_{x \in V, \| x\| =1} \| T x \| \\
 &  \leq & \max_{x \in V_1, \| x\| =1} \| S x \| +
\max_{x \in V_2, \| x\| =1} \| T x \|\\
&=& s_p(S) + s_q(T).
\eea
It follows from Proposition \ref{min max prop} that for $r = \mbox{codim } V$,
$s_{r+1}(S+T) \leq s_p(S) + s_q(T)$.
As $r +1\leq p+q -1$, $s_{p+q -1}(S+T) \leq s_{r+1}(S+T)$.
$\Box$

\begin{prop}\label{read sing values}
Let $H_1, H_2$ be two Hilbert spaces.
Let $\{u_n: n \geq 1 \}$ be a Hilbert basis of $H_1$ 
and $\{v_n: n \geq 1 \}$ be a Hilbert basis of $H_2$.
Let $F$ be the subspace of $H_1$ spanned by finite linear 
combinations of the $u_n$, $n\geq 1$ and $a_n, n \geq 1$ 
a sequence of real numbers. 
Define a linear operator $T: F \ri H_2$ by setting
$T u_n = a_n v_n$. Assume that $a_n$ is decreasing and converges to zero. Then $T$ can be uniquely extended to a linear operator
from $H_1$ to $H_2$. $T$ is compact and $s_n(T) = a_n$.
\end{prop}
\textbf{Proof:}
It is clear that $F$ is dense in $H_1$ and $T$ is uniformly continuous on $F$, thus a unique extension to $H_1$ exists.
As $a_n$ converges to zero, it can be shown that $T$ is the strong
 limit of a sequence of finite dimensional range linear operators $T_k $ defined
by $ T_k u_n = a_n v_n$, if $n \leq k$, $T_k u_n =0$,
if $n > k$. 
It can be checked that $T^* v_n = a_n u_n $
thus  $T^* T u_n = a_n^2 u_n$ and $s_n(T) = a_n$.
$\Box$

\section{Neumann eigenvalues and inverse Laplacian}
\subsection{Weyl's asymptotic formula}
Let $\Omega$ be a bounded open subset of $\RR^d$
and $0 < \lambda_1 \leq \lambda_2 \leq ...$
be the sequence of the non-zero eigenvalues for the Neumann problem 
in $\Omega$.
Weyl \cite{weyl1912asymptotische} showed the asymptotic formula
\bean \label{weyl}
\lambda_n \sim 4 \pi^2 (\omega_d |\Omega|)^{-\f{2}{d}} 
n^{\, \f{2}{d}},
 \eean
where $\omega_d$ is the volume 
of the unit ball in $\RR^d$ and 
$|\Omega|$ is the Lebesgue measure of $\Omega$.
Some regularity on the boundary $\p \Omega$ is required 
for \eqref{weyl} to hold.
In fact, it is known \cite{hempel1991essential}
that for any closed subset $[0,\infty)$ and
$d \geq 2$,
there is a
bounded open set of $\RR^d$,  such that the essential spectrum
of the Neumann Dirichlet operator is equal to that closed set.
Courant's proof of Weyl's formula \cite{CourantHilbert}
involves covering 
$\Omega$ by a countable family of cubes and  was originally written
for a domain with a $C^2$ boundary.
Netrusov and Safarov 
provide in the introduction
of their paper on Weyl's formula and domain regularity
\cite{netrusov2005weyl}
a concise explanation for
why Courant's proof applies to general Lipschitz domains, even
in the Neumann case.
In summary, formula \eqref{weyl} is valid if $\Omega$ 
is a Lipschitz domain and we will assume in the rest of this paper 
that $\Omega$ is a connected domain that is at least Lipschitz regular.

\subsection{Neumann eigenvectors and associated compact operators}
\label{f_nn mu_n}
Let $L^{2,0}(\Omega)$ be the subspace of functions in $L^2(\Omega)$ that are orthogonal to constants.
Introduce the Sobolev subspace $H^{1,0}(\Omega) =
H^1(\Omega) \cap L^{2,0}(\Omega)$.
Due to our assumption on the geometry of $\Omega$,
the Neumann eigenvalues for the Laplace operator satisfy
\eqref{weyl}.
Setting $\mu_n = \lambda_n^{-1}$, the sequence $\mu_n$
is valued in $(0,\infty)$, is decreasing, and converges to zero.
Let $f_n$ be an eigenfunction in $H^{1,0}(\Omega)$
such that
\bean
\Delta f_n + \mu_n^{-1} f_n =0 \label{neu1}\\
\nabla f_n \cdot \nu =0 \label{neu2}\\
\int_{\Omega} f_n^2 =1 \label{neu3}.
\eean 
Due to Poincar\'e's inequality we can choose the inner product
on $H^{1,0}(\Omega)$ to be defined by 
$<f,g> = \int_\Omega \nabla f \cdot \nabla g$.
Note that 
\bean
\int_{\Omega} \nabla f_n \cdot \nabla \varphi = 
\mu_n^{-1} \int_{\Omega}  f_n  \varphi, \label{weak neu}
\eean
for all $\varphi $ in $H^{1}(\Omega)$ since \eqref{weak neu}
also holds for constant functions $\varphi$.
We can choose $f_n$ such that the set $\{ f_n: n \geq 1 \}$
forms a Hilbert basis of $L^{2,0}(\Omega)$. 
Since the functions $f_n$ are orthonormal in $L^{2,0}(\Omega)$,
it follows from \eqref{weak neu} that the 
functions $\mu_n^{\f12} f_n$ 
are orthonormal in $H^{1,0}(\Omega)$. 
It can actually be proved that  
$\{ \mu_n^{\f12} f_n: n \geq 1 \}$
forms a Hilbert basis of $H^{1,0}(\Omega)$.\\
We now define an inverse Laplacian.
Set 
\bea
J: L^{2,0}(\Omega) \ri L^{2,0}(\Omega),\\
 J f_n  = \mu_n f_n.
\eea
$J$ defines a compact operator thanks to 
Proposition \ref{read sing values}.
Note that in our argument this operator will play the same role
as the operator $J$ defined in 
\cite[p. 416]{ha1986} or \cite[Chap III, p. 120]{GohbergKrein}.
 The operator $J$ is in some sense 
an inverse Laplacian since the range of 
$J$ is in $H^{1,0}(\Omega)$ and the rule
$\int_{\Omega} \nabla J f_n \cdot \nabla \varphi = 
\int_{\Omega}  f_n \varphi$ for all $\varphi$ in  $H^{1,0}(\Omega)$
extends to 
\bean \label{inv Lapla}
\int_{\Omega}   \nabla J f \cdot  \nabla \varphi= 
\int_{\Omega}  f \varphi ,
\eean
 for all $f$ in $L^{2,0}(\Omega)$ and
$\varphi$ in  $H^{1}(\Omega)$.\\
Next, using Proposition \ref{read sing values}
we define the compact operator
\bea
T: L^{2,0}(\Omega) \ri H^{1,0}(\Omega),\\
 T f_n  = \mu_n f_n.
\eea
According to Proposition \ref{read sing values}, 
$s_n(T) =  \mu_n^{\f12} $.
We then introduce for $1 \leq k \leq d$  a   derivative type operator
\bea
D_k: H^{1,0}(\Omega) \ri L^{2,0}(\Omega),  \\
 D_k f_n  = \p_k f_n - \f{1}{|\Omega|} \int_\Omega \p_k f_n .
\eea
Note that
\bea
\int_\Omega | D_k f_n|^2 \leq \int_\Omega  |\p_k f_n |^2
\eea
since $D_k f_n$ is the orthogonal projection of $\p_k f_n$
on $L^{2,0} (\Omega)$. 
From there, it is clear that $D_k$ is continuous and $\| D_k \| \leq 1$.
Finally we set 
\bean \label{Nkdef}
N_k = D_k T.
\eean 
Thanks to Proposition \ref{norm and sing} we can write
the estimate  for $1 \leq k \leq d$,
\bean
s_n(N_k) \leq \mu_n^{\f12}. \label{s_n(N_k)} 
\eean

From \eqref{inv Lapla}, we can now write the integration by parts formula
\bea
\int_{\Omega}   \varphi f =
\int_{\Omega} \nabla \varphi \cdot  \nabla J f = 
\sum_{i_1=1}^d  \int_{\Omega} \p_{i_1} \varphi N_{i_1} f,
\eea
 for all $f$ in $L^{2,0}(\Omega)$ and
$\varphi$ in  $H^{1}(\Omega)$.\\
By iteration, if $p \geq 1$ is an integer and $\varphi $
is in $H^{r}(\Omega)$,
\bean \label{int by parts}
\int_{\Omega}   \varphi f =
\sum_{1 \leq i_1, ..., i_p \leq d} \int_{\Omega} \p_{i_1} ... \p_{i_p}
\varphi N_{i_p} ... N_{i_1} f,
\eean
for $f$ in $L^{2,0}(\Omega)$. 

\subsection{Higher norms of Neumann eigenvectors}
If $\p \Omega$ is $C^2$ regular, then $f_n$ defined by 
(\ref{neu1}-\ref{neu3}) is $H^2$ Sobolev regular
and $\| f_n \|_{H^2(\Omega)} \leq C \mu_n^{-1}$, where
$C$ only depends on $\Omega$. This is due to the fact 
that elliptic regularity all the way to the boundary is also valid for Neumann problems
\cite[Chapter~IV]{Mikhailov1978}.\\
If $p \geq 2$ is an integer and $\p \Omega$ is $C^p$ regular,
 there is a constant $C$ depending only on $\Omega$ and $p$
such that if $\Delta u = f$ in $\Omega$ and $\nabla f \cdot \nu =0$ on 
$\p \Omega$ then
\bean
 \| u \|_{H^p(\Omega)} \leq C \| f \|_{H^{p-2}(\Omega)}.
\eean
From there    we can show by induction using 
(\ref{neu1}-\ref{neu3})
that
\bean \label{r reg}
\| f_n \|_{H^p(\Omega)} \leq C \mu_n^{-\f{p}{2}}.
\eean
Next, let $a_n$ be a sequence of real numbers and $0<q_1<q_2$
two integers.
As 
\bea
\Delta (\sum_{n=q_1}^{q_2} a_n f_n) = \sum_{n=q_1}^{q_2} a_n \mu_n^{-1}
f_n,
\eea
it follows from \eqref{r reg}, (\ref{neu1}-\ref{neu3}),
and the orthogonality of the functions $f_n$  
that
\bea
\| \sum_{n=q_1}^{q_2} a_n f_n \|_{H^2(\Omega)}^2 \leq 
C^2 \sum_{n=q_1}^{q_2} a_n^2 \mu_n^{-2}
\eea
If $p\geq 3$, using that the gradients $\nabla f_n$  are orthogonal in
 $L^2$, we similarly find that
\bea
\| \sum_{n=q_1}^{q_2} a_n f_n \|_{H^3(\Omega)}^2 \leq 
C^2 \sum_{n=q_1}^{q_2} a_n^2 \mu_n^{-3}.
\eea
More generally, a simple argument by induction shows that 
\bean \label{in Hr norm}
\| \sum_{n=q_1}^{q_2} a_n f_n \|_{H^p(\Omega)}^2 \leq 
C^2 \sum_{n=q_1}^{q_2} a_n^2 \mu_n^{-p}.
\eean

\section{Proof of Theorems stated in section \ref{intro} }
\subsection{Singular values of integral operators with $L^2$ kernel. 
Translations, dilations, and set inclusions.}
\begin{thm} \label{only L2}
Let $\Omega$ be an open set of $\RR^d$ and
$ \Omega'$ be an open set of $\RR^{d'}$.
Let $k$ be in $L^2(\Omega' \times \Omega) $.
Define the linear operator
\bean \label{as in }
K:  L^2(\Omega) \ri L^2(\Omega') , \no \\
K \varphi(x) = \int_\Omega k(x,y) \varphi(y) dy.
\eean
Then $K$ is compact and 
\bean \label{sing2}
\sum_{n \geq 1} s_n(K)^2 = \int_{\Omega' \times \Omega}
 k(x,y)^2 dx dy
\eean
\end{thm}
\textbf{Proof:}
This theorem is  proved in textbooks in the case where 
$\Omega = \Omega'$ and $K$ is symmetric non-negative. 
We cover a more general case for  sake of completeness.
Let $\varphi$ be in $L^2(\Omega)$.
As
\bea
\int_{\Omega}  (\int_{\Omega'} |k(x,y) \varphi (y)| d y)^2 dx \leq
(\int_{\Omega'} \varphi^2 (y) dy)  (\int_{\Omega' \times \Omega}
 k(x,y)^2 dx dy ),
\eea
	$K \varphi$ is in $L^2(\Omega')$ and $K$ is a continuous
	linear operator from $L^2(\Omega)$ to $L^2(\Omega')$.
	To show that $K$ is compact we note that the argument
	in section II.14 of \cite{gohberg2013basic}
	which was made in the case where $\Omega=\Omega'$, 
	can be trivially generalized to the present case. Indeed,
	let $\{\varphi_m: m\geq 1 \}$ be a Hilbert basis 
	of $L^2(\Omega)$ and $\{\psi_n: n\geq 1 \}$ be a Hilbert basis 
	of $L^2(\Omega')$. Then 
	$\{\varphi_m \psi_n: m\geq 1, n\geq 1 \}$ is a Hilbert basis 
	of $L^2(\Omega' \times \Omega)$.
	Defining 
	\bea
	k_N (x,y) = \sum_{1 \leq m, n \leq N}
	 (\int_{\Omega' \times \Omega}     k(x',y')  \psi_n(x') \varphi_m(y') dx' dy')
	\psi_n(x) \varphi_m(y) 
	\eea
	and
	\bea
K_N:  L^2(\Omega) \ri L^2(\Omega') , \\
K_N\varphi(x) = \int_\Omega k_N(x,y) \varphi(y) dy,
\eea
$K_N$ has finite dimensional range. 
Since $\|k_N - k \|_{L^2(\Omega' \times \Omega)}$
converges to zero as $N \ri \infty$, 
 $K_N$ converges in operator norm to $K$,
thus $K$ is compact. \\ 
Let $u_n$ be orthonormal functions in $L^2(\Omega)$
such that $K^*K u_n = s_n(K) u_n$. 
If $s_n(K)\neq 0$, define
 $v_n = K u_n/s_n(K)$. The functions $v_n$ form an
 orthonormal set in $L^2(\Omega')$.  It follows
that
\bean \label{on tensor basis}
\delta_{m,n} s_n(K) = 
\int_{\Omega' \times \Omega} k(x,y) v_m(x) u_n(y)
dx dy.
\eean
If the range of $K$ is not dense, we  complete  the set of functions $v_n$ to obtain a Hilbert basis of $L^2(\Omega')$.
We can then prove that $\{  v_m(x) u_n(y): n \geq 1, m  \geq 1 \}$
is a Hilbert basis of $L^2(\Omega' \times \Omega)$. 
As $k$ is in  $L^2(\Omega' \times \Omega)$, formula 
\eqref{sing2} follows from \eqref{on tensor basis}.
$\Box$\\
	Under the conditions of Theorem \ref{only L2}, using Proposition
	\ref{prop:decreasing} proved in appendix it follows that
	\bea
	\lim_{ n \ri \infty} n^{\f12} s_n(K) = 0.
	\eea

	\begin{prop} \label{dila trans}
	Let $k$ and $K$ be defined as in Theorem \ref{only L2}.
	Let $t>0$ be  a scalar,  $v_0$ be in $\RR^d$,
	and $\tilde{\Omega} = t \Omega + v_0$.
	Set $\tilde{k} (x,v)= k(x, \f{v - v_0}{t} )$
	and $\tilde{K}$ the integral operator $L^2(\tilde{\Omega})
	\ri L^2(\Omega')$  defined from $\tilde{k} $. The singular values 
	of $K$ and $\tilde{K}$ 
	are related by the formula
	\bean \label{dila trans form}
	s_n(\tilde{K} ) = t^{\f{d}{2}} s_n(K).
	\eean
	\end{prop}
	\textbf{Proof:}
	Let $\varphi_n$ be an eigenvector for the eigenvalue 
	$s_n(K)^2$  of $K^*K$.
	As
	\bea
	s_n(K)^2\varphi_n(y) =  K^*K \varphi_n(y) = 
	\int_{\Omega'} k(x,y) \int_{\Omega} k(x,u) \varphi_n(u) du  dx,
	\eea
	setting $u = \f{v - v_0}{t}$, $y = \f{w- v_0}{t}$ we find,
	\bea
	t^{d} s_n(K)^2\varphi_n(\f{w- v_0}{t}) =  
	\int_{\Omega'} k(x,\f{w- v_0}{t}) \int_{\tilde{\Omega}} k(x, \f{v - v_0}{t}) \varphi_n( \f{v - v_0}{t})  dv  dx.
	\eea
	Since this formula holds for any arbitrary $\Omega$, $t>0$, and $v_0 \in
	\RR^d$, formula \eqref{dila trans form} is proved.
	$\Box$\\

	\begin{prop} \label{inclusion}
	Let $k$ and $K$ be defined as in Theorem \ref{only L2}.
	Let $\Omega_2$ be an open subset of $\Omega$ such that
	$\Omega_2 \subset \Omega$ and define
	\bea 
K_2:  L^2(\Omega_2) \ri L^2(\Omega') , \no \\
K_2 \varphi(x) = \int_{\Omega_2} k(x,y) \varphi(y) dy.
\eea
	Then
	\bean  \label{inclusion form} 
	s_n(K_2) \leq s_n(K).
	\eean
		\end{prop}
	\textbf{Proof:}
	Let $\varphi_1, ..., \varphi_{n-1} $ be orthonormal
	in $L^2(\Omega)$
	 and such that $K^*K \varphi_j = s_j(K)^2 \varphi_j$, $j=1, ..., n-1$.
	Let
	\bea
	X= \{ \psi \in L^2(\Omega):
	<\psi,\varphi_1>_{L^2(\Omega)} = ... = 
	<\psi,\varphi_{n-1}>_{L^2(\Omega)} = 0 \}.
	\eea
	Then by Proposition \ref{min max prop}
	\bea
	s_n(K) = \max_{\psi \in X, \|  \psi  \|_{L^2(\Omega)} =1}
	\| K \psi \|.
	\eea
	We note that the co-dimension
	\bea
	X_2= \{ \psi \in L^2(\Omega_2):
	<\psi,\varphi_1>_{L^2(\Omega_2)} = ... = 
	<\psi,\varphi_{n-1}>_{L^2(\Omega_2)} = 0 \}
	\eea
	is less or equal than $n-1$.
	Moreover, any element in $X_2$ with norm 1 can be extended by zero to become an element
	of $X$ with norm 1.
	Applying again  Proposition \ref{min max prop}
	\bea
	s_n(K_2) \leq  \max_{\psi \in X_2, \| \psi  \|_{L^2(\Omega_2)} =1}
	\| K_2 \psi \| \leq  \max_{\psi \in X, | \psi  \|_{L^2(\Omega)} =1}
	\| K \psi \| = s_n(K) .
	\eea
	$\Box$\\

	\subsection{Proof of Theorem \ref{thm1}}
Let $k$ be in $L^2(\Omega') \times H^p(\Omega)$ as in the statement
of Theorem \ref{thm1}. 
Let $f$ be in $L^2_0(\Omega)$. 
Note that the co-dimension of $L^2_0(\Omega)$ in $L^2(\Omega)$
is 1.
According to  \eqref{int by parts},
\bean \label{bigsum}
\int_{\Omega}   k(x,y) f(y) dy =
\sum_{1 \leq i_1, ..., i_p \leq d} \int_{\Omega} \p_{i_1} ... \p_{i_p}
 k(x,y) N_{i_p} ... N_{i_1} f(y) dy,
\eean
	where all the derivatives are in the $y$ variable.
	We then set 
	\bea
K_{i_1, ..., i_p}:  L^2_0(\Omega) \ri L^2(\Omega') , \\
K_{i_1, ..., i_p} \varphi(y) = \int_\Omega\p_{i_1} ... \p_{i_p}
 k(x,y)  \varphi(y) dy.
\eea
	Let $K_0$ be the restriction of $K$
	to $ L^2_0(\Omega)$. 
	We can write the  operator identity
	\bea
	K_{0 } = \sum_{1 \leq i_1, ..., i_p \leq d} K_{i_1, ..., i_p} N_{i_p} ... N_{i_1}.
	\eea
	Given \eqref{pq}, \eqref{s_n(N_k)}, and the simplification
		\bea
		(p-1)(n-1)+ 1 + n -1 = p(n-1)  +1, 
		\eea
		we can prove by induction on $p$
	that for $n \geq 1$
	that $s_{p(n-1) + 1 } ( N_{i_p} ... N_{i_1})\leq \mu_n^{\, \f{p}{2}}$
	thus 
	\bean \label{multip N}
	s_{(p+1)(n-1) + 1 } ( K_{i_1, ..., i_p} N_{i_p} ... N_{i_1})\leq \mu_n^{\, \f{p}{2}}
	s_n(K_{i_1, ..., i_p}).
	\eean
	Note that the sum in \eqref{bigsum} has $d^p$ terms.
	Let $L_j$, $1 \leq j \leq d^p$, be compact operators.
	We can prove 
	by induction using \eqref{pq plus}
	that
	\bean \label{big sum s}
	s_{d^p(n-1) + 1 } (\sum_{j=1}^{d^p} L_j)
	\leq \sum_{j=1}^{d^p} s_n(L_j).
	\eean
	Combining \eqref{multip N} and \eqref{big sum s}
	we find that
	\bea
	s_{d^p  (p+1) (n-1)+1 } (K_0) \leq \mu_n^{\, \f{p}{2}} 
	\sum_{1 \leq i_1, ..., i_p \leq d}
	s_n (K_{i_1, ..., i_p} ).
	\eea
	As the co-dimension of $L^2_0(\Omega)$ in $L^2(\Omega)$ is 1, by 
	Proposition \ref{restriction}, this implies that
	\bean \label{interm}
	s_{d^p  (p+1) (n-1) + 2} (K) \leq \mu_n^{\, \f{p}{2}} \sum_{1 \leq i_1, ..., i_p \leq d}
	s_n (K_{i_1, ..., i_p} ).
	\eean
	Recall that $\mu_n = \lambda_n^{-1} = C_1 n^{- \f{2}{d}}$
	where $C_1>0$ depends only on $d$ and on $|\Omega|$ and is given by
	\eqref{weyl}. Set $C_2 = \sum_{1 \leq i_1, ..., i_p \leq d}
	s_n (K_{i_1, ..., i_p} ) $ and $C_3= d^p (p+1)$.
	Recalling that singular values are taken in decreasing order we obtain
	from \eqref{interm}
	\bea
	s_{C_3 n} (K) \leq C_1^{\, \f{p}{2}} C_2 n^{-\f{p}{d}}. 
	\eea
	For $m > C_3$, we write $ m = C_3 n + \alpha$, where
	$0 \leq \alpha < C_3$.
	Then
	\bea
	 s_m(K) \leq  C_1^{\, \f{p}{2}}  C_3^{\, \f{p}{d}} (m- C_3 +1)^{-\f{p}{d}} \sum_{1 \leq i_1, ..., i_p \leq d}
	s_n (K_{i_1, ..., i_p} ) , \quad C_3= d^p (p+1), \quad n=[ \f{m}{C_3} ],
	\eea
	and estimate  \eqref{main est} is proved.\\
	To infer estimate \eqref{particular 1} from \eqref{main est}
	we recall that $d$ and $p$ are fixed.
	For a fixed integer $m \geq 1 $ the set 
	$ \{ n \in \NN: m =[\f{n}{d^p(p+1)}] \}$
	contains $d^p(p+1)$ elements.
	Next, the sum $\ds \sum_{1 \leq i_1, ..., i_p \leq d}
	s_m (K_{i_1, ..., i_p} )$ has $d^p$
	terms and   
$\sum_{m \geq 1} 	s_m (K_{i_1, ..., i_p} )^2 < \infty$
since  $\p_{i_1} ... \p_{i_p}
 k$ is in $L^2(\Omega' \times \Omega)$, see Theorem 
\ref{only L2}. \\
Finally, \eqref{particular 2} results from \eqref{particular 1} 
thanks to Proposition \ref{prop:decreasing} in appendix.
$\Box$

	
	\subsection{Proof of Theorem \ref{thm2} }\label{proof of thm2}
	The first part  of Theorem \ref{thm2} pertains to  Lipschitz regular
	open sets $\Omega$.
	By Proposition \ref{dila trans} we may assume without loss
	of generality that $\Omega_2 = (0,\pi)^d \subset \Omega$.
	We set 
	\bean \label{phi_m_j def}
	\varphi_{m_j} (y_j)= \s{\f{2}{\pi}}\cos(m_j y_j),
	\eean
where $m_j \geq 1$ is an integer.
Note that
the functions 
\bea
\vp_{m_1}(y_1) ...  \vp_{m_d}(y_d), \quad m_1 \geq 1, ... ,  m_d \geq d
\eea
form an orthonormal system in $L^2(\Omega_2)$.
If we fix the integers $j_1 \geq 0, ..., j_d \geq 0 $
 	the functions 
\bean  \label{also ortho}
m_1^{-j_1}... m_d^{-j_d}
\vp_{m_1}^{(j_1)}(y_1) ...  \vp_{m_d}^{(j_d)}(y_d), \quad m_1 \geq 1, ... ,  m_d \geq d
\eean 
also form an orthonormal system in $L^2(\Omega_2)$, where
$\vp_{m_i}^{(j_i)}$ is the derivative of $\vp_{m_i}$
of order $j_i$.
For $m$ in $\NN^d$ we use the notation
	$m=(m_1, ..., m_d)$, $|m| = m_1+ ... +m_d$.
For $(x,y)$ in $\Omega \times \Omega '$ and $\lambda$ in
$\RR$ define
\bea
k(x,y) = \sum_{n \geq 1} \sum_{|m|=n} 
n^{-p +\lambda} \vp_{m_1}(y_1) ...  \vp_{m_d}(y_d)
f_{m_1, ..., m_d}( x),
\eea
where the functions $f_{m_1, ..., m_d}$
are orthonormal in $L^2(\Omega')$ and the parameter 
$\lambda$ is such that
\bean \label{conv}
 \sum_{n \geq 1} \sum_{|m|=n} 
n^{2 \lambda} < \infty.
\eean 
Since
	\bea 
	\mbox{card } \{ m \in \NN^d: |m| = n \} = 
	\f{(n+d-1) ... (n+1)}{(d-1)!} 
	\eea
	it suffices that $\lambda < - \f{d}{2}$
	for \eqref{conv} to hold.
Due to the orthonormal property of the system \eqref{also ortho}
and the requirement  \eqref{conv} it is clear
that $k$ is in $L^2(\Omega')\times H^p(\Omega_2)$.
Due to symmetry and periodicity properties of the functions
$\vp_{m_i}$, $k$ is also in $L^2(\Omega')\times H^p(\Omega)$.
Now, define $k_2$ and $K_2$ from $\Omega$ and $\Omega_2$
as in Proposition \ref{inclusion}.
Given that   
 	the functions 
  \eqref{also ortho} form an orthonormal set, due to Proposition 
\ref{read sing values} the set of  singular values of $K_2$
is given by $ 
n^{-p +\lambda}$ with the repetition pattern
 $|m|=n $, $m \in \NN^d$. In other words, the singular value
$n^{-p +\lambda}$ is repeated $\f{(n+d-1) ... (n+1)}{(d-1)!} $ times.
An argument by induction on $n$ can show that 
	\bean \label{counting}
	 \sum_{  0 \leq j \leq n } \f{(d+j-1)!}{j! (d-1)!} = \f{(n+d)!}{n! d!}=
	\f{(n+d) ... (n+1)}{d!}.
	\eean
	As $ \f{n^{d}}{d!} \leq \f{(n+d) ... (n+1)}{d!}  $,
it follows that  
\bea
s_{[\f{n^{d}}{d!}]} (K_2)\geq n^{-p + \lambda}
\eea
Let $m$ be an integer such that $\f{(n-1)^d}{d!} \leq m < \f{n^d}{d!}$.
Then $s_m(K_2) \geq n^{-p+\lambda}$,
thus
\bea
s_m(K_2) \geq (m^{\f{1}{d}} (d!)^{\f{1}{d}} + 1)^{-p+\lambda}
\eea
and by Proposition \ref{inclusion},
\bea
s_m(K) \geq (m^{\f{1}{d}} (d!)^{\f{1}{d}} + 1)^{-p+\lambda}.
\eea
Estimate \eqref{main est 2.1} is proved by taking
$\lambda$ arbitrarily closed to $-\f{d}{2}$.\\
To prove the second part Theorem \ref{thm2} ,
	let  $p \geq 2$ be  an integer and assume  that $\p \Omega$ is $C^p$ regular.
	Let $\{ \varphi_n: n\geq 1 \}$ be a Hilbert basis of $L^2(\Omega')$.
	Let $f_n$ and $\mu_n$ be as in section
	\ref{f_nn mu_n}.
	Let $a_n>0$ be such that
		\bean
		\mbox{
	$  a_n$ is
	decreasing,} \label{decreas}\\
\sum_{n \geq 1} \mu_n^{-p} a_n^2  < \infty, \label{series cond}\\
\mbox{and for some }C_a>0, \f{a_n}{a_{2n}} \leq C_a, \mbox{ for all } n\geq 1.
\label{not too fast}
	\eean
	For example, $a_n= \mu_n^{\, \f{p}{2}} \f{1}{n} $ 
	satisfies  (\ref{decreas}-\ref{not too fast}).
 Define the function in $L^2(\Omega \times \Omega')$
\bean \label{k def}
k(x,y) = \sum_{n \geq 1} a_n \varphi_n(x) f_n(y),
\eean
	and the operator $K:  L^2(\Omega) \ri L^2(\Omega')$, 
	$
K f(x) = \int_\Omega k(x,y) f(y) dy.$
By Proposition \ref{read sing values}, $s_n (K) = a_n$.
	\begin{lem} \label{series lem}
	The kernel $k$ defined by \eqref{k def} is in
	$L^2(\Omega') \times H^p(\Omega)$ and for $1 \leq i_1,... ,i_p \leq d$
	\bean
	\p_{i_1} ... \p_{i_p} k(x,y) =  
	\sum_{n \geq 1} a_n \varphi_n(x) \p_{i_1} ... \p_{i_p}f_n(y),
	\label{der sum}
	\eean
	where all derivatives are in the $y $ variable and the sum 
	in \eqref{der sum}
	converges 
	in $L^2(\Omega' \times \Omega)$.
	In particular the series
	\bean \label{L2 series}
	\sum_{n \geq 1} a_n  \p_{i_1} ... \p_{i_p}f_n
	\eean
	converges in $L^2(\Omega)$.
	\end{lem}
	\textbf{Proof:}
	Let  $0<q_1<q_2$ be
two integers.
For a fixed $x$, using \eqref{in Hr norm}, 
\bea 
\int_\Omega |\sum_{q_1\leq n \leq q_2}
a_n  \p_{i_1} ... \p_{i_p}f_n(y) \varphi_n(x)  dy|^2
dy \leq C^2 \sum_{q_1\leq n \leq q_2} a_n^2 \mu_n^{-p} \varphi_n(x)^2,
\eea
thus,
\bea 
\int_{\Omega'} \int_\Omega |\sum_{q_1\leq n \leq q_2}
a_n  \p_{i_1} ... \p_{i_p}f_n(y) \varphi_n(x)  dy|^2 dx
dy \leq C^2 \sum_{q_1\leq n \leq q_2} a_n^2 \mu_n^{-p} ,
\eea
and the series 
	$\sum_{n \geq 1} a_n \varphi_n(x) f_n(y)$
	converges in $L^2(\Omega') \times H^p(\Omega)$
	thanks to the 
	Cauchy criterion and \eqref{series cond}. 
	Formula \eqref{der sum} follows in the sense of weak derivatives.
	$\Box$\\\\
	
	We now show that the decay rate estimate 
	\eqref{main est} is optimal for the integral operator $K$ 
	with  kernel $k$ defined in \eqref{k def}.
	We note that from Lemma \ref{series lem}
	$K_{i_1, ..., i_p}^* \varphi_n = a_n \p_{i_1} ... \p_{i_p}f_n$.
	Let $X=\{ \varphi_1, ..., \varphi_{n-1} \}^\perp$ in $L^2(\Omega')$.
	By \eqref{dual} and \eqref{min max},
	\bea
	s_n(K_{i_1, ..., i_p}) = 	s_n(K_{i_1, ..., i_p}^*)  \leq
\max_{\psi \in X, \| \psi \| =1} \| K_{i_1, ..., i_p}^* \psi \|.
	\eea
	For $\psi $ in $X$, set $\psi = \sum_{k \geq n} b_k \varphi_k$ with 
	 $\sum_{k \geq n} b_k^2 =1$. 
	Then \bea
	K_{i_1, ..., i_p}^* \psi = \sum_{k \geq n } b_k a_k \p_{i_1} ... \p_{i_p}f_k
	\eea
	and using \eqref{in Hr norm}
	\bea
	\| K_{i_1, ..., i_p}^* \psi \|_{L^2(\Omega)}^2 \leq C^2
	 \sum_{k \geq n } b_k^2  a_k^2 \mu_k^{-p}.
	\eea
	Since $a_k^2 \mu_k^{-p}$ is positive and decreasing we find that,
	\bea
	\max_{\psi \in X, \| \psi \| =1} \| K_{i_1, ..., i_p}^* \psi \| \leq 
	C a_n \mu_n^{-\f{p}{2}},
	\eea
	thus 
	\bean
	s_n(K_{i_1, ..., i_p}) \leq C a_n \mu_n^{-\f{p}{2}}.
	\label{s_n(Ki_1, ..., i_p) }
	\eean
	From \eqref{not too fast} it follows by induction that for any integer 
	$r \geq 1$, 
	\bea
	\f{a_n}{a_{2^r n}} \leq C_a^r.
	\eea
	As previously, we denote $C_3=d^p(p+1)$.
	We fix $r$ such that $C_3 +1  \leq 2^r$, so by \eqref{decreas}
	\bean
	\f{a_n}{a_{(C_3+1) n}} \leq C_a^r. \label{final Ca}
	\eean
	By \eqref{s_n(Ki_1, ..., i_p) }
	\bea
	s_{[\f{n}{C_3}]}(K_{i_1, ..., i_p}) \leq C a_{[\f{n}{C_3}]}
	\mu_{[\f{n}{C_3}]}^{-\f{p}{2}},
	\eea
	but
	\bea
	 a_{[\f{n}{C_3}]} \leq C_a^r 
	a_{(C_3+1)[\f{n}{C_3}]} \leq  C_a^r  a_n.
	\eea
	for all $n$ large enough.
	Similarly,  since  $\mu_n = \lambda_n^{-1}$,
	recalling \eqref{weyl}, 
	$	\mu_{[\f{n}{C_3}]}^{-\f{p}{2}}$ is bounded 
	above by a constant independent of $n$ times 
	$\mu_n^{-\f{p}{2}}$.
	Altogether,
	we find that
	\bean
	s_{[\f{n}{C_3}]}(K_{i_1, ..., i_p}) \leq C_{a, \mu} a_n \mu_n^{-\f{p}{2}},
	\label{s_n(Ki_1, ..., i_p) 2}
	\eean
	where $C_{a, \mu} >0$ does not depend on $n$.
	Recalling $a_n = s_n(K)$, \eqref{weyl}, and $\mu_n = \lambda_n^{-1}$,
	summing over all $i_1, ..., i_p$ in $\{1,...,d\}$
 formula \eqref{main est 2}	is proved. 
The second part of Theorem \ref{thm2} is now proved. $\Box$

	\subsection{The real analytic case: proof of the first part of Theorem 
	\ref{thm3}} \label{The real analytic case}
	We first precisely define what it means for
	 $k(x,y)$ to be in $L^2(\Omega' \times \Omega)$
	and real analytic in $y$ in a neighborhood of $\ov{\Omega}$.
	For $y=(y_1, ..., y_d)$ in $\RR^d$ and $I=(i_1, ..., i_d)$
	in $\NN^d$ denote the scalar
	\bea
	y^I= y_1^{i_1} ... y_d^{i_d}.
	\eea
	We also use the notation $|y|_\infty = \max \{ |y_1|, ..., |y_d| \}$.
	We say that 
	 $k(x,y)$ is in $L^2(\Omega' \times \Omega)$
	and real analytic in $y$ in a neighborhood of $\ov{\Omega}$
	if there is an open set $U$ of $\RR^d$ such that 
	$\ov{\Omega} \subset U$ and for all $y'$ in $U$ 
	there is a multi-indexed sequence $a_{I,y'}(x)$ in $L^2(\Omega')$
	and a radius$R_{y'} >0$
	such that 
	\bea
	\sum_{I \in \NN^d} \| a_{I,y'}\|_{L^2(\Omega')} < \infty,
	\eea
	and
	\bea
	 k(x,y) = \sum_{I \in \NN^d}  a_{I,y'}(x) (y-y')^I,
	\eea
	for all $y$ in $\RR^d$ such that $|y - y'|_\infty < R_{y'}$.
	For example, this occurs if $k(x,y)$ is real analytic for $(x,y)$ 
	in a neighborhood of  $\ov{\Omega' \times \Omega}$.\\
	Since $\ov{\Omega}$ is compact we can 
	use a partition of unity $\sum_{j=1}^q \eta_j = 1 $
on  $\ov{\Omega}$ where each $ \eta_ j  \geq 0$ is smooth and supported
in one of the balls with center $y'$ and radius $R_{y'}/3$.
Let $\by_j$ and $R_{\by_j}$ denote this center and this radius. 
	Accordingly, we set 
	\bean
	k_j(x,y) = k(x,y) \eta_j(y). \label{kj def}
	\eean
	Set for $z$ in $\CC^d$
	\bea
	M_j = \max_{|z- \by_j| = \f23 R_j} \| k_j(x,z ) \|_{L^2(\Omega')}. 
	\eea
	If $z$ in $\CC^d$ is such that $|z - \by_j |_\infty < \f23 R_j$,
	\bea
	k(x,z ) =
	\f{1}{(2i\pi)^d} \int_{{\cal C}_d } ... \int_{{\cal C}_1}
	\f{k(x,\zeta_1, ..., \zeta_d)}{ (\zeta -z )^{(1, ..., 1)}},
	\eea
	where the circle ${\cal C}_k $ in the complex plane is
	centered at the $k$-th coordinate of $\by_j$ and its radius is
	$\f23 R_j$.
	It follows that 
	\bea
	a_{I, \by_{j}} (x) =
	\f{1}{(2i\pi)^d} \int_{{\cal C}_d } ... \int_{{\cal C}_1}
	\f{k_j(x,\zeta_1, ..., \zeta_d)}{ (\zeta - \by_j )^{I+(1, ..., 1)}},
	\eea
	from where we derive the Cauchy estimate
	\bean \label{cauchy}
	\| a_{I, \by_{j}}\|_{L^2(\Omega')} \leq
	M_j (\f23 R_j)^{- |I|},
 	\eean
	where $|I| = i_1 + ... + i_d$.\\
		We now estimate the singular values of $k_j$ generalizing the  idea 
	used in the one-dimensional  case ($d=d'=1$) 
	\cite[Theorem 15.20]{kress1989linear}, \cite{little1984eigenvalues}.
	Define the partial sum
	\bea
	k_j^N(x,y) =\sum_{  I \in \NN^d, |I|\leq N} a_{I, \by_{j}} ( y - \by_j)^I.
		\eea
	It follows from \eqref{cauchy} that if $| y - \by_j | \leq \f13 R_j$,
	\bea
	\| k_j(x,y) - k_j^N(x,y)\|_{L^2(\Omega')} 
	&\leq& M_j \sum_{ |I|\geq N+1} 2^{-|I|}\\
		&=& M_j \sum_{r \geq N+1} \f{(r+d-1) ... (j+1)}{(d-1)!} 2^{-r} \\
	&\leq & M_j  (\f32)^{-N},
	\eea
	for all $N$ large enough.
	Recall that if $| y - \by_j | \geq \f13 R_j$, then 
	$ k_j(x,y) = k_j^N(x,y)=0$.
	Since the first singular value is equal to the operator norm,
	\bean \label{s1 est}
	s_1( k_j  - k_j^N) 
	\leq M_j |\Omega|^{\f12} (\f32)^{-N}.
	\eean
	The operator $k_j^N$ has finite dimensional range. 
	Let $\alpha_N$ be the dimension of that range.
	Given that the dimension of the space of homogeneous polynomials in $d$ variables
	with degree $i$ is $\f{(d+i-1)!}{i! (d-1)!}$,
	an argument by induction on $n$ can show that 
	\bea
	 \sum_{  0 \leq i \leq N } \f{(d+i-1)!}{i! (d-1)!} = \f{(N+d)!}{N! d!},
	\eea
	so $\alpha_N$ is bounded above by a polynomial in $N$ with degree less or equal to 
	$d$.
	Let $K_j$ be the integral operator associated to $k_j$
	and $K_j^N$ be the integral operator associated to $k_j^N$ 
	as in \eqref{as in }.
	Note that $K_j^N$ has at most $\alpha_N$ non-zero singular values.
	Then by \eqref{pq} and \eqref{s1 est}, 
	\bea
	s_{\alpha_N + 1} (K_j)\leq s_1 (K_j - K_j^N) +
		s_{\alpha_N + 1} (K_j^N)
		\leq  M_j |\Omega|^{\f12} (\f32)^{-N}
	\eea
	For all $N$ large enough $\alpha_N \leq N^d$,
	thus $s_{N^d } (K_j)\leq C (\f32)^{-N}$, for a constant $C$ independent of $N$.
	It follows that for any integer $m$ such that 
	$N^d \leq m < (N+1)^d$,
	\bea
	s_m(K_j) \leq C (\f32)^{-N} \leq \f32 C (\f32)^{-m^{\f{1}{d}}}
	\eea
	By  \eqref{pq plus} and \eqref{kj def}
	\bea
	s_{qm}(K) \leq \f32 qC (\f32)^{-m^{\f{1}{d}}}
	\eea
	Now for integers $n$ such that $qm \leq n< q(m+1)$,
	given that  $- m^{\f{1}{d}} \leq - (m+1)^{\f{1}{d}} +1$,
	\bea
	s_n(K) \leq s_{qm}(K) \leq \f94 q C (\f32)^{-(m+1)^{\f{1}{d}}},
	\eea
	 it follows that 
	\bea
	s_n(K) \leq s_{qm}(K) \leq \f94 q C \tau^{n^{\f{1}{d}}},
	\eea
	where $\tau =(\f32)^{ -q^{-\f{1}{d} } } $ is in $(0,1)$. This proves the first part
	of Theorem \ref{thm3}.
	
	\subsection{Proof of the second part
	of Theorem \ref{thm3}}

	\subsubsection{The  case where 
	 $\Omega$ is the rectangle  $(0,\pi)^d$}
	Here too, we use the functions
	$\varphi_{m_j}$ defined in 
	 \eqref{phi_m_j def}.
	Accordingly,  we label
	\bea
	\bigcup_{j=1}^d \bigcup_{m\geq 1} \{f_m^j\}
	\eea
	a Hilbert basis of $L^2(\Omega')$.
	Fix $\tau$ in $(0,1)$. 
	We then set 
	\bean \label{k anal}
	k(x,y) & =&
	(\sum_{m_1 \geq 1} \tau^{m_1} \varphi_{m_1}(y_1) f_{m_1}^1(x))
	...
	(\sum_{m_d \geq 1} \tau^{m_d} \varphi_{m_d}(y_d) f_{m_d}^d(x))\\
	 & =& \sum_{n\geq 1} \sum_{|m| =n }
	\tau^n  \varphi_{m_1}(y_1) ...  \varphi_{m_d}(y_d) f_{m_1}^1(x) ...
	f_{m_d}^d(x). \no
	\eean
	It is clear that $k$ is real analytic in $y$, valued in $L^2(\Omega')$
	since $\tau \in (0,1)$ implies that 
	$\sum_{m_j \geq 1} \tau^{m_j} \varphi_{m_1}(y_j) $
	is holomorphic for $y_j$ in a strip
	centered on the real axis.
	Let $K$ be the operator $K:  L^2(-1,1)^d \ri L^2(\Omega')$
	defined as in \eqref{as in }.
	Since
	\bea
	\mbox{card } \{ m \in \NN^d: |m| = n \} = 
	\f{(n+d-1) ... (n+1)}{(d-1)!} \geq [\f{n^d}{(d-1)!}],
	\eea
	It follows that
	\bea
	s_{[\f{n^d}{(d-1)!}]} (K)\geq \tau^{n}.
	\eea
	Now let $m$ be an integer such that
	$\f{(n-1)^d}{(d-1)!}  \leq m < \f{n^d}{(d-1)!} $.
	Then $s_m(K) \geq \tau^{  n}$, so
	\bean \label{on a square}
	s_m(K) \geq \tau (\tau^{((d-1)!)^{\f{1}{d}}})^{m^{\f{1}{d}}},
	\eean
	Since $\tau \in (0,1)$, $\tau^{((d-1)!)^{\f{1}{d}}}$
	can be arbitrarily close to 1 and the statement of Theorem
	\ref{thm3} is proved in this case.
	
	\subsubsection{Case where $\Omega$ is an arbitrary open set}
	By proposition \ref{dila trans} 
	we may assume that 
$(0,\pi)^d\subset \Omega$. 
Then $k$ defined by \eqref{k anal} is analytic in 
$\RR^d$.  We use  Proposition \ref{inclusion} and we consider
 inequality \eqref{inclusion form} 
where $k$ and $K$ are related  as in Theorem \ref{only L2},
$\Omega_2 = (0,\pi)^d$ and $K_2$ is defined from $K$
as in  Proposition \ref{inclusion}. 
Recalling estimate \eqref{on a square},
the proof of Theorem \ref{thm3} is now complete.
	$\Box$

	\subsection{Manifolds}
Suppose ${\cal M} $ is a compact $d$-dimensional
  manifold in $\RR^{\tilde{d}}$ of class $C^p$
	or an open and bounded manifold in
	$\RR^{\tilde{d}}$ of class $C^p$
	with a $C^p$ boundary.
Let 	 $d \sigma$ denote the surface measure on ${\cal M} $.
Let 	 ${\cal M}' $ be a  manifold in $\RR^{d'}$.
	For $(x,y)\ri k(x,y)$ in $L^2({\cal M} ') \times H^p({\cal M} )$,
define the linear operator
\bea 
K:  L^2(\Omega) \ri L^2(\Omega') , \no \\
K \varphi(x) = \int_{{\cal M}} k(x,y) \varphi(y) d\sigma (y).
\eea
Using  $C^p$ diffeormorphisms locally mapping  
${\cal M}$ to open sets in $\RR^d$ and using associated partitions of unity
such as those defined in section \ref{The real analytic case},
it is  clear  that
	Theorems \ref{thm1} and \ref{thm2} 
	can be extended to the case of manifolds. 
	
\appendix
\section{Appendix: decreasing sequences, convergent series, and convergence rates}
Proposition \ref{prop:decreasing} proved in this appendix appears in 
\cite[page 122]{gohberg1969}, at least in the case where $q=1$.
We provide here an alternative proof which we believe to be simpler. In 
particular, it only uses elementary results on series.
\begin{lem} \label{integer lemma}
  Let $p>0$ be a real number.  Let $b_k$ be a 
	strictly increasing sequence in $\mathbb{N}$.
	Then
    \[ \sum_{k=1}^\infty \left( 1 - \left( \frac{b_k}{b_{k+1}} \right)^p \right) \]
    diverges.
\end{lem}

\textbf{Proof:}
    Arguing by contradiction, assume that
    $ \sum_{k=1}^\infty \left( 1 - \left( \frac{b_k}{b_{k+1}} \right)^p \right) $
    converges. Then\\
    $ \lim_{k \to \infty} \left( 1 - \left( \frac{b_k}{b_{k+1}} \right)^p \right) = 0. $
    Set $a_k = 1 - \left( \frac{b_k}{b_{k+1}} \right)^p$. 
  It follows that  $\lim_{k \to \infty} \frac{-\ln(1-a_k)}{a_k} = 1$, so 
		 $\sum_{k=1}^\infty -\ln(1-a_k)$ converges.
		Note that 
    \bea -\ln(1-a_k) = 
		-\ln \left( \left( \frac{b_{k}}{b_{k+1}} \right)^p \right) = 
		p\ln(b_{k+1}) - p\ln(b_k),
		\eea
    thus 
    \[ \sum_{k=1}^n -\ln(1-a_k)  = \sum_{k=1}^n \left( p\ln(b_{k+1}) - p\ln(b_k) \right) = p\ln(b_{n+1}) - p\ln(b_1). \]
    But $b_{n+1} \geq n +1$
    so as $n \ri \infty$ this  contradicts that $\sum_{k=1}^n -\ln(1-a_k) $ converges. 
$\Box$


\begin{prop}\label{prop:decreasing}
    Let $p \geq 0$ and $q>0$ be two real numbers. Let $a_n > 0$ be a decreasing sequence such that
    \[ \sum_{n=1}^\infty n^p a_n^q \]
    converges. Then
    \[ \lim_{n \to \infty} n^{\f{p+1}{q}} a_n = 0. \]
\end{prop}

\textbf{Proof:}
    Arguing by contradiction, suppose
  $ n^{\f{p+1}{q}} a_n $ does not converge to zero. 
    Then there exists $\epsilon > 0$ and a subsequence $a_{n_k}$ such that
    $  n_k^{\f{p+1}{q}} a_{n_k} \geq \epsilon $, that is,
		\bean \label{ank}
		 a_{n_k}^q \geq \frac{\epsilon^q}{n_k^{p+1}},
		\eean
    for all $k \geq 1$. 
		
		Therefore
    \bea
		\sum_{n_k  < j \leq n_{k+1} } j^p a_j^q 
		\geq  \sum_{n_k < j \leq n_{k+1}  } j^p a_{n_{k+1}}^q 
	\geq  \frac{\epsilon^q }{n_{k+1}^{p+1}}
	\sum_{n_k < j \leq n_{k+1} }  j^p.
		\eea
		But
		\bea
		\sum_{n_k < j \leq n_{k+1} }  j^p \geq
		\int_{n_k}^{n_{k+1}} x^p dx  = \f{1}{p+1} ( n_{k+1}^{p+1} - n_{k}^{p+1}),
		\eea
		thus
		\bea
		\sum_{n_k  < j \leq n_{k+1} } j^p a_j^q 
		\geq \f{\epsilon^q}{p+1} (
		 1 - ( \frac{n_k}{n_{k+1}} )^{p+1} ) .
		\eea
		Summing for $k \geq 1$ and using   Lemma \ref{integer lemma}
		this contradicts that $\sum_{n=1}^\infty n^p a_n^q$
		converges. $\Box$



\end{document}